\documentclass{article}
\usepackage{amsmath}
\usepackage{amsfonts}
\usepackage{amssymb}
\usepackage{graphicx}

\def\qed{$\,_{\square}$}

\newtheorem{remark}{Remark}
\newtheorem{theo}{Theorem}
\newtheorem{lemma}{Lemma}
\newtheorem{conj}{Conjecture}

\begin{document}
\title{Recursions and tightness for the maximum of the discrete, two
dimensional Gaussian Free Field}
\author{Erwin Bolthausen\thanks{Work supported in part by 
SNSF project 200020-125247 and 
the Humboldt foundation.}\,, Universit\"{a}t Z\"{u}rich
\and Jean Dominique Deuschel\thanks{Work supported in part by 
DFG-Forschergruppe 718.}\,, Technische Universit\"{a}t Berlin
\and Ofer Zeitouni\thanks{Work supported in part by NSF grant
DMS-0804133.}\,, University of Minnesota and Weizmann Institute}
\date{June 25, 2010}
\maketitle
\begin{abstract}
	We consider the maximum of the discrete two dimensional 
	Gaussian free field in a box, and prove the existence of
	a (dense) deterministic subsequence along which 
	the maximum, centered at its mean, is tight; this still leaves open
	the conjecture that tightness holds without the need for subsequences.
	The method of proof relies on an argument developed by Dekking and Host
	for branching random
	walks with bounded increments and on comparison results 
	specific to Gaussian fields.
\end{abstract}
\section{Introduction and main result}
We consider the discrete Gaussian Free Field (GFF)  in a two-dimensional
box of side $N+1$, with Dirichlet boundary conditions. That is, let
$V_N=([0,N]\cap \mathbb{Z})^2$, $V_N^o=((0,N)\cap \mathbb{Z})^2$,
and let
$\{w_m\}_{m\geq 0}$ denote a simple random walk started 
in $V_N$ and killed at $\tau=\min\{m: w_m\in
\partial V_N\}$ (that is, killed upon hitting the boundary
$\partial V_N=V_N\setminus V_N^o$). 
For $x,y\in V_N$, define
$G_N(x,y)=E^x(\sum_{m=0}^\tau {\bf 1}_{w_m=y})$, where $E^x$ denotes expectation
with respect to the random walk started at $x$. The GFF is
the zero-mean Gaussian field $\{X_z^N\}_z$
indexed by $z\in V_N$ with covariance
$G_N$.

Let $X_N^*=\max_{z\in V_N} X_z^N$.
It was proved in \cite{BDG} that $X_N^*/(\log N)\to c$ with $c=2
\sqrt{2/\pi}$, and the proof is closely related to
the proof of the law of large numbers for the
maximal displacement of a branching random walk (in $\mathbb{R}$). 
Based on the analogy with the maximum of independent
Gaussian variables and the case of branching random walks, the following
is a natural conjecture.
\begin{conj}
	\label{conj-1}
The sequence of random variables $Y_N:=X_N^*-EX_N^*$ is tight.
\end{conj}
To the best of our knowledge, the sharpest result in this direction is due
to \cite{sourav}, who shows that the variance of $Y_N$ is $o(\log N)$;
in the same paper, Chatterjee also analyzes related Gaussian fields,
but in all these examples, does not prove tightness.
We defer to Section \ref{sec-bibref} for some pointers to the relevant
literature concerning the Gaussian free field 
and the origin of Conjecture \ref{conj-1}.

The goal of this note is to prove a weak form of the conjecture. Namely,
we will prove the following.
\begin{theo}
\label{theo-main}
There is a deterministic sequence $\{N_k\}_{k\geq 1}$ such that
the sequence of random variables $\{Y_{N_k}\}_{k\geq 1}$ is tight.
\end{theo}
More information on the sequence $\{N_k\}_{k\geq 1}$ is provided 
below in Remark \ref{remm-1}.

It is of course natural to try to improve the tightness from subsequences
to the full sequence. As will be clear from the proof, for that it is enough
to prove the existence of a constant $C$ such that $EX_{2N}^*\leq EX_N^*+C$.
This is weaker than, and implied by, the conjectured behavior of $EX_N^*$, 
which is
\begin{equation}
	\label{eq-top}
	EX_N^*=c\log N - c_2 \log \log N+O(1),
\end{equation}
for $c=2\sqrt{2/\pi}$ and an appropriate  $c_2$,
see e.g. \cite{CLD} and Remark \ref{rem-2}.

Finally, although we deal here exclusively with the GFF, it should be clear
from the proof that the analysis applies to a much wider class of models.

\section{Preliminary considerations}
Our approach is motivated by the proof of tightness of branching random
walks (BRW) with independent increments, in the spirit of \cite{DH}
(see also the argument in \cite{BZ}). We will thus first 
introduce a branching-like structure in the GFF. Unfortunately,
this structure is not directly suitable for analysis, and so we 
later modify it.

\subsection{The basic branching structure} 
We consider $N=2^n$ in what follows, write $Z_n=X_N^*$ 
and identify an integer
$m=\sum_{\ell=0}^{n-1}
m_{i+1} 2^i$ with 
its binary expansion $(m_{n},m_{n-1},\ldots,m_1)$. 
For $k\geq 1$,
introduce 
the sets of $k$-diadic integers
$$A_k=\{m\in \{1,\ldots,N\}: m=(2l+1)N/2^k \; \quad 
\mbox{\rm for some integer $l$}\}.$$
Note that if $m\in A_k$ then $m_i=0$ for $i\leq n-k$ and $m_{n-k}=1$.
Then, define the $\sigma$-algebras
$${\cal A}_k=\sigma(X_z^N: z=(x,y), x \, \mbox{\rm or}\,
y\in \cup_{i\leq k} A_i).$$
Finally, for every   $z=(x,y)\in V_N^o$, 
write $z_i=(x_i,y_i)$ with $x_i,y_i$ 
denoting as above the $i$th digit in the binary expansion of $x,y$.
We introduce the random variables
\begin{equation}
\label{eq-1}
\xi^{z_1,\ldots,z_k}_{z_{k+1},\ldots,z_n}= E[X_z^N|{\cal A}_k].
\end{equation}
We then have the decomposition
\begin{equation}
\label{eq-2b}
X_z^N=\xi^{z_1}_{z_2,\ldots,z_n}+X^{z_1}_{z_2,\ldots,z_n},
\end{equation}
where, by the Markov property of the GFF, the collections
$\{X^{z_1}_{\cdot}\}_{z_1\in V_1}$
are i.i.d.  copies of the GFF in the box $V_{N/2}$, and are independent
of the collection of random variables
$\{\xi^{z_1}_{\cdot}\}$.
Iterating, we have the representation
\begin{equation}
\label{eq-2}
X_z^N=\xi^{z_1}_{z_2,\ldots,z_n}+
\xi^{z_1,z_2}_{z_2,\ldots,z_n}+\ldots+
\xi^{z_1, z_2,\ldots,z_{n-1}}_{z_n}\,,
\end{equation}
where all the summands in the right side of \eqref{eq-2} are independent.
Recall that $X_N^*=\max_{z\in V_N} X_z$.

We can now explain the relation with branching random walks: should the 
random variables in the right side of \eqref{eq-2} not depend on the 
conditioning (that is, the superscript), \eqref{eq-2} would correspond
precisely to a branching random walk (on a four-ary tree), with time-dependent 
increments. For such 
BRW, a functional recursion for the law
of $X_N^*$ can be written down, and used 
to prove tightness (see \cite{BZ1} and \cite{ABR}).
Unfortunately, no such simple functional recursions are available
in the case \eqref{eq-2}. For this reason, 
we first modify the representation \eqref{eq-2}, and
then adapt an argument of \cite{DH}, originally
presented in the context of BRW.
To explain our goal, note that
we have from \eqref{eq-2b} that
\begin{equation}
\label{eq-3}
Z_n=X_N^*=\max_{z\in V_1}((X_{N/2}^*)^z+D^{z,N}),
\end{equation}
where the variables $\{(X_{N/2}^*)^z\}_z$ are four i.i.d. copies
of $X_{N/2}^*$, and the $D^{z,N}$ are complicated fields but
$D^{z,N}_{z_2,\ldots,z_n}\geq \min_{z_2,\ldots,z_n} \xi^z_{z_2,\ldots,z_n}$.
Unfortunately,
the $D^{z,N}$ variables are
far from being uniformly 
bounded, and this fact prevents the application of the
argument from \cite{DH}.

\subsection{Two basic lemmas}
In this subsection we present two preliminary lemmas that will allow 
for a comparison of the GFF between different scales. The first
shows that the maximum of the sum of two  zero mean fields tends to 
be larger than each of the fields.

\begin{lemma}
\label{lem-new1}
Let $\{X_i\}_{i\in V_Ni}$ and  $\{Y_i\}_{i\in V_N}$ be two independent
random fields and assume that $EY_i=0$. Then,
\begin{equation}
\label{eq-new1}
E\max_{i\in V_N} (X_i+Y_i)\geq 
E\max_{i\in V_N} X_i\,.
\end{equation}
\end{lemma}
\noindent
{\bf Proof} Let $\alpha\in V_N$ be such that
$\max_{i\in V_N} X_i =X_\alpha$ (in case several $\alpha$s satisfy the above
equality, choose the first according to lexicographic order).
We then have
$$E\max_{i\in V_N} (X_i+Y_i)\geq 
E(X_\alpha+Y_\alpha)=EX_\alpha+EY_\alpha=
EX_\alpha=
E\max_{i\in V_N} X_i\,,$$
where the second equality is due to the independence of the
fields and the fact that $EY_i=0$ for all $i$. \qed

By \eqref{eq-3} and Lemma \ref{lem-new1},
we have that 
\begin{equation}
\label{eq-wis1}
EZ_{n+1}\geq EZ_n\,.
\end{equation}
The following lemma gives 
a control in the opposite direction.
\begin{lemma}
        \label{lem-oldeq3}
        There exists a sequence $n_k\to \infty$ and a constant $C$ such that
        \begin{equation}
                \label{eq-4}
                EZ_{n_k+1}\leq EZ_{n_k}+C\,.
        \end{equation}
\end{lemma}
{\bf Proof} From \cite{BDG} there exists a constant $c>0$ so that
$EZ_n/n\to c$.
Fixing arbitrary
$K$  and defining
$I_{n,K}=\{i\in [n,2n]: EZ_{n+1}>EZ_n+K\}$,
one has from \eqref{eq-wis1}
and the existence of the limit
$EZ_n/n\to c$ that
$$ \limsup_{n\to\infty}\frac{|I_{n,K}|}{2n}
\leq \frac{c}{K}\,.$$
In particular, choosing $K=3c$ it follows that for all $n$ large, there
exists an $n'\in [n,2n]$ so that
$$EZ_{n'}\leq EZ_{n'-1}+K\,,$$
as claimed.\qed
 
\subsection{Proof of Theorem \ref{theo-main}
}
By \eqref{eq-3} and Lemma \ref{lem-new1}, we get that
$$EX_{2N}^*\geq E\max_{z\in V_1}(X_N^*)^z\geq E\max(X_N^*,\tilde X_N^*)\,,$$
where $\tilde X_N^*$ is an independent copy of $X_N^*$.
Using the equality
$\max(a,b)=(a+b+|a-b|)/2$, we get that
$$ EZ_{n+1}-EZ_n\geq E|X_N^*-\tilde X_N^*|/2\,.$$
For the sequence $n_k$ and the constant $C$ from 
Lemma \ref{lem-oldeq3}, we thus
get that
$$2C\geq E|X_{2^{n_k}}^*-E\tilde X_{2^{n_k}}^*|\,.$$
This shows that the sequence
$\{X_{2^{n_k}}^*\}_{k\geq 1}$ is tight and completes the proof of Theorem
\ref{theo-main}.
\qed

\begin{remark}
\label{remm-1}
The subsequence $n_k$ provided in 
Lemma \ref{lem-oldeq3} can be taken with
density arbitrary close to $1$, as can be seen from the following modification
of the proof. 
Fixing arbitrary
$K$ and $\epsilon$ and defining
$I_{n,\epsilon,K}=\{i\in [n,n(1+\epsilon)]: EZ_{n+1}>EZ_n+K\}$,
one has from \eqref{eq-wis1}
 and the existence of the limit
$EZ_n/n\to c$ (with $c=2\sqrt{2/\pi}$) that 
$$ \limsup_{n\to\infty}\frac{|I_{n,\epsilon,K}|}{n\epsilon}
\leq \frac{c}{K}\,.$$
It is of course
of interest to see whether one can take $n_k=k$.
Minor modifications of the proof of Theorem \ref{theo-main}
would then yield Conjecture
\ref{conj-1}.
\end{remark}
%
\begin{remark}
Minor modifications of the proof of Theorem \ref{theo-main}
also show that
if there exists a constant $C$ so that $EX^*_{2N}\leq EX^*_{N}+C$ for all
integer $N$, then Conjecture \ref{conj-1} holds. 
\end{remark}
\begin{remark}
	\label{rem-2}
	For Branching Random Walks, under suitable assumptions
	it was established in \cite{ABR} that
	\eqref{eq-top} holds.
	Running the argument above then immediately 
	implies the tightness of the minimal (maximal) displacement, centered
	around its mean.
\end{remark}
	
\section{Some bibliographical remarks}
\label{sec-bibref}
The Gaussian free field has been extensively
studied in recent years, in both its continuous and discrete forms. For
an accessible review, we refer to \cite{Sheff}. The fact that the GFF 
has a logarithmic decay of correlation invites a comparison with 
branching random walks, and through this analogy a form of
Conjecture \ref{conj-1} is implicit in \cite{CLD}. 
This conjecture is certainly
``folklore'', see e.g.  open problem \#4 in \cite{sourav}. 
For some one-dimensional
models (with logarithmic decay of correlation) 
where the structure of the maxima can be analysed, we refer to 
\cite{BF,LDF}. The analogy with branching random walks has been
reinforced by the study of the so called {\it thick points} of the GFF, both
in the discrete form \cite{daviaud} and in the continuous form \cite{miller}.


{\bf Acknowledgment} 
A proof of Theorem \ref{theo-main} was first provided in \cite{BDZ}.
In that proof some additional estimates controlling the
field $D^{z,N}$, see \eqref{eq-3}, are provided.	
Yuval Peres observed that that these estimates are redundant, and suggested
the proof presented here. We thank Yuval for this observation and for his 
permission to use it here.


\end{document}